\documentclass[a4paper, 10pt]{amsart}

\usepackage{etoolbox}
\makeatletter
\let\ams@starttoc\@starttoc
\makeatother
\usepackage[parfill]{parskip}
\makeatletter
\let\@starttoc\ams@starttoc
\patchcmd{\@starttoc}{\makeatletter}{\makeatletter\parskip\z@}{}{}
\makeatother

\usepackage{color,soul}
\definecolor{red}{rgb}{1,0,0}
\setulcolor{red}

\usepackage{amssymb}
\usepackage{amsmath}
\usepackage{amsthm}
\usepackage{fancyhdr}
\usepackage{esint}
\usepackage{enumerate}

\usepackage[colorlinks=true,linkcolor=blue,citecolor=blue,urlcolor=blue]{hyper ref}
\usepackage{bookmark}

\swapnumbers
\newtheorem{lemma}{Lemma}[section]
\newtheorem{prop}[lemma]{Proposition}
\newtheorem{thm}[lemma]{Theorem}
\newtheorem{cor}[lemma]{Corollary}

\theoremstyle{definition}
\newtheorem{defn}[lemma]{Definition}
\newtheorem{example}[lemma]{Example}
\newtheorem{rem}[lemma]{Remark}
\newtheorem{ass}[lemma]{Assumption}

\numberwithin{equation}{section}

\renewcommand{\(}{\left(}
\renewcommand{\)}{\right)}
\renewcommand{\~}{\tilde}

\newcommand{\bs}{\backslash}

\newcommand{\cn}{\colon}

\newcommand{\N}{\mathbb{N}}
\newcommand{\R}{\mathbb{R}}

\renewcommand{\a}{\alpha}
\renewcommand{\b}{\beta}
\newcommand{\g}{\gamma}
\renewcommand{\d}{\delta}
\newcommand{\e}{\epsilon}
\renewcommand{\k}{\kappa}
\renewcommand{\l}{\lambda}
\renewcommand{\O}{\Omega}

\newcommand{\D}{\Delta}

\newcommand{\p}{\varphi}

\newcommand{\fr}[2]{\frac{#1}{#2}}

\newcommand{\Thm}{\begin{thm}}
\newcommand{\eThm}{\end{thm}}
\newcommand{\Def}{\begin{defn}}
\newcommand{\eDef}{\end{defn}}
\newcommand{\Prop}{\begin{prop}}
\newcommand{\eProp}{\end{prop}}
\newcommand{\Rem}{\begin{rem}}
\newcommand{\eRem}{\end{rem}}
\newcommand{\Lem}{\begin{lemma}}
\newcommand{\eLem}{\end{lemma}}
\newcommand{\eq}{\begin{equation}}
\newcommand{\eeq}{\end{equation}}
\newcommand{\Ex}{\begin{example}}
\newcommand{\eEx}{\end{example}}
\newcommand{\pf}{\begin{proof}}
\newcommand{\epf}{\end{proof}}
\newcommand{\Cor}{\begin{cor}}
\newcommand{\eCor}{\end{cor}}
\newcommand{\Ass}{\begin{ass}}
\newcommand{\eAss}{\end{ass}}
\newcommand{\SAl}{\begin{align}\begin{split}}

\newcommand{\hra}{\hookrightarrow}

\newcommand{\mrm}{\mathrm}

\newcommand{\hp}{\hphantom}

\protected\def\ignorethis#1\endignorethis{}
\let\endignorethis\relax

\begin{document}

\title[Quantitative oscillation estimates]{Quantitative oscillation estimates for almost-umbilical closed hypersurfaces in Euclidean space}
\author{Julian Scheuer}
\date{\today}
\thanks{This work is being supported by the DFG}
\address{Dr. Julian Scheuer\ \\ Ruprecht-Karls-Universit\"at, Institut f\"ur Angewandte Mathematik, Im Neuenheimer Feld 294, 69120 Heidelberg, Germany}
\email{scheuer@math.uni-heidelberg.de}
\begin{abstract}
We prove $\e$-closeness of hypersurfaces to a sphere in Euclidean space under the assumption that the traceless second fundamental form is $\d$-small compared to the mean curvature. We give the explicit dependence of $\d$ on $\e$ within the class of uniformly convex hypersurfaces with bounded volume. 
\end{abstract}
\maketitle
\tableofcontents

\section{Introduction}
In this article we investigate the potential of the traceless second fundamental form, also called the \textit{umbilicity tensor}, 
\eq \mathring{A}=A-\fr{\mrm{tr}(A)}{n}g\eeq of a hypersurface embedded in the Euclidean space
to pinch other geometric quantities of the hypersurface. Questions like this arise from the well known fact, that $\mathring{A}=0$ implies that the hypersurface must be a sphere. Then it is natural ask if this behavior is a kind of continuous, in the sense that a small traceless second fundamental form implies closeness to a sphere. During the last decade, substantial progress has been made towards a better understanding of this question. In 2005 an article by Camillo de Lellis and Stefan M{\"u}ller  appeared,
\cite{DeLellisAndMueller2005}, where the estimate 
\eq \inf_{\l\in\R}\|A-\l g\|_{L^{2}(M)}\leq C\|\mathring{A}\|_{L^{2}(M)}\eeq
was proven for hypersurfaces $M\subset\R^{3}.$ From this, the authors deduced $W^{2,2}$-closeness to a sphere. One year later, in \cite{DeLellisAndMueller2006}, the authors made a step towards uniform closeness and showed that in addition the metric is $C^{0}$-close to the standard sphere metric. In 2011, one of de Lellis' PhD students, Daniel Perez, proved in the class of hypersurfaces with volume 1 and bounded second fundamental form, that for given $\e>0$ there exists $\d>0,$ such that a $\d$-small traceless second fundamental form yields $\e$-closeness to a sphere, compare \cite[Cor. 1.2]{PerezDiss}. He used an argument via contradiction and it does not seem possible to extract the $\e$-dependence of $\d$ along his proof. In \cite[p. xvi]{PerezDiss} the author posed the derivation of a quantitative dependence as an open problem. In this article we tackle this problem and prove the following theorem.

\Thm\label{Main}
Let $n\geq 2$ and $X\cn M^{n}\hra\R^{n+1}$ be the smooth, isometric embedding of a closed, connected, orientable and mean-convex hypersurface. Let $0<\a<1.$ Then there exists $c>0,$ such that whenever we have $\e<c|M|^{\fr{1}{n}}$ and the pointwise estimate 
\eq \label{Pinch} \|\mathring{A}\|\leq H|M|^{-\fr{2+a}{n}}\e^{2+\a}\eeq
holds, then $M$ is strictly convex and 
\eq M\subset B_{\sqrt{\fr{n}{\l_{1}(M)}}+\e}(x_{0})\bs B_{\sqrt{\fr{n}{\l_{1}(M)}}-\e}(x_{0}).\eeq
The constant $c$ depends on $n,$ $\a,$ $\|\~{A}\|_{\infty}$ and $\|\~{A}^{-1}\|_{\infty},$ where $|M|=\mrm{vol}(M),$ $\~{A}=|M|^{\fr{1}{n}}A,$ $\l_{1}(M)$ is the first nonzero eigenvalue of the Laplace-Beltrami operator on $M$ and $x_{0}$ is the center of mass of $M.$
\eThm
 
A more detailed description of the notation involved here will be presented in section \ref{Notation}. Thus in the class of uniformly convex hypersurfaces of unit volume we obtain $\e$-closeness to a sphere, if $\mathring{A}$ is of order $\e^{2+\a}$ and $\e$ is sufficiently small.

The author's motivation to find a quantitative dependence like this arose from his work on inverse curvature flows in the Euclidean space. In \cite[Appendix A]{Schnuerer02/2006} Oliver Schn{\"u}rer derived a pinching estimate of the traceless second fundamental form for hypersurfaces evolving by the inverse Gauss curvature flow in $\R^{3}.$ Ben Andrews applied estimates like this to bound the difference between circumradius $r_{+}$ and inradius $r_{-}$ of the surface in \cite[Section 4]{Andrews1999}. However, we are not aware whether those methods may be transferred to higher dimensions. Clearly, Theorem \ref{Main} provides an estimate of $r_{+}-r_{-}$ in terms of $\mathring{A}.$ Indeed, we are going to apply this estimate to prove asymptotical roundness of hypersurfaces solving an inverse curvature flow equation in $\R^{n+1},$ cf. \cite{js:Ipcf}.

Let us give an overview over the main ingredients involved in the proof. Certainly we need a result, which somehow yields the transition from \textit{qualitative to quantitative}.
We found the following result due to Julien Roth. We formulate a special case and only the statements which are of interest to our proof. 

\Thm \cite[Thm. 1]{Roth2007}\label{Roth}
Let $(M^{n},g)$ be a compact, connected and oriented Riemannian manifold without boundary isometrically immersed in $\R^{n+1}.$ Assume that $|M|=1$ and $H_{2}>0.$ Then for any $p\geq 2$ and $\e>0$ there exists a constant $C_{\e}=C_{\e}(n,\|H\|_{\infty},\|H_{2}\|_{2p}),$ such that if
\eq  \label{PC} \l_{1}(M)\(\int_{M}H\)^{2}-n\|H_{2}\|_{2p}^{2}>-C_{\e}\eeq
is satisfied, then 
\eq M\subset B_{\sqrt{\fr{n}{\l_{1}}}+\e}(x_{0})\bs B_{\sqrt{\fr{n}{\l_{1}}}-\e}(x_{0}),\eeq
where $x_{0}$ is the center of mass of $M$  and $H_{2}$ is the second normalized elementary symmetric polynomial.
\eThm

This theorem  is a generalization of \cite{ColboisAndGrosjean09/2006} to higher $k$-th mean curvatures. There are generalizations to ambient spaces of bounded sectional curvature, cf. \cite{GrosjeanAndRoth2012}, as well.
At first glance, it does not seem to be a quantitative result, but a rather tedious scanning of the proof shows, that $C_{\e}$ can be chosen to be of order $\e^{2},$ compare section \ref{Qualitative closeness revisited}. 

Certainly, this $\e^{2}$ gives insight into the question, where the order $\e^{2+\a}$ comes from in Theorem \ref{Main}. It is an interesting question, whether, and if how, this could be improved.

Thus we have to derive (\ref{PC}) from (\ref{Pinch}). Firstly, we need to relate the first eigenvalue of the Laplacian to the traceless second fundamental form. This transition has another stop at the Ricci tensor. The following result, due to Erwann Aubry, relates the Ricci tensor to $\l_{1}.$ It was proven in \cite{Aubry2007}, but is accessible more easily in \cite[Thm. 1.6]{Aubry2009}. Again, we only cite the aspects, which are relevant to our work.

\Thm \cite[Thm. 1.6]{Aubry2009}\label{Aubry}
For any $p>\fr{n}{2}$ there exists $C(n,p),$ such that if $M^{n}$ is a complete manifold with
\eq \label{Aubry1}\int_{M}(\underline{\mrm{Ric}}-(n-1))_{-}^{p}<\fr{|M|}{C(n,p)},\eeq
then $M$ is compact and satisfies 
\eq \l_{1}(M)\geq n\(1-C\(\fr{1}{|M|}\int_{M}(\underline{\mrm{Ric}}-(n-1))_{-}^{p}\)^{\fr 1p}\).\eeq
Here, $\underline{\mrm{Ric}}=\underline{\mrm{Ric}}(x)$ denotes the smallest eigenvalue of the Ricci tensor at $x\in M$ and for $y\in\R$ we set $y_{-}=\max(0,-y).$ 
\eThm

The other quantities in (\ref{PC}) can be controlled with the help of (\ref{Pinch}) quite easily. Thus the only ingredient, which is left, is to control the Ricci tensor in (\ref{Aubry1}). The following result, due to Daniel Perez, \cite{PerezDiss} and also to De Lellis and M{\"u}ller for $n=2,$ \cite{DeLellisAndMueller2005}, is helpful.

\Thm \cite[Thm. 1.1]{PerezDiss}\label{Perez}
Let $n\geq 2,$ $p\in(n,\infty)$ and $c_{0}>0.$ Then there exists $C(n,p,c_{0})>0,$ such that for any smooth, closed and connected hypersurface $M\subset\R^{n+1}$ with
\eq |M|=1\eeq
and
\eq \|A\|_{p}\leq c_{0}\eeq
we have 
\eq\label{Perez1} \min_{\mu\in\R}\|A-\mu g\|_{p}\leq C\|\mathring{A}\|_{p}.\eeq
\eThm

This result will enable to move, via the Ricci tensor, to an estimate on $\l_{1}$ and to finally provide the estimate (\ref{PC}). Then the result follows. There largest technical difficulty is, that we finally need $L^{\infty}$ bounds, where the theorems \ref{Aubry} and \ref{Perez} only make statements on $L^{p}$ norms. We will present the way to handle this in section \ref{Quantitative spherical closeness}.

Note, that we will not need to know the explicit value of $\mu_{0}$ in (\ref{Perez1}), where the minimum is attained. However, this is another interesting question with some history. According to \cite[p. 50]{PerezDiss}, it was Gerhard Huisken to suggest an inverse mean curvature flow approach to prove, that the minimum is attained at 
\eq \mu=\fr{1}{|M|}\int_{M}H.\eeq 
In \cite[p. 52, Ch. 3.4]{PerezDiss} this is proven for $n\geq 2,$ $p=2$ and for closed convex hypersurfaces. Unfortunately, the case $p=2$ is not enough in our case. Hence we have to deal with the little technical difficulty, that $\mu_{0}$ is not explicitly known. 

Now we start with the detailed analysis of the problem at hand and start with an explanation of our notation.

\section{Notation and preliminaries}\label{Notation}

In this article we consider closed embedded hypersurfaces $M^{n}\subset\R^{n+1}.$ We follow the notation as it appears in the references as closely as possible.

 $g=(g_{ij})$ denotes the induced metric of $M^{n},$ $A=(h_{ij})$ the second fundamental form and $\k_{i},$ $i=1,\dots,n,$ the principal curvatures ordered pointwise,
\eq \k_{1}\leq \dots \leq \k_{n}.\eeq
The volume of $M$ is
\eq |M|=\int_{M}1\ d\mu,\eeq
where $\mu$ is the canonical surface measure associated to $g.$

 $\l_{1}(M)$ denotes the first nonzero eigenvalue of $-\D,$ where $\D$ is the Laplace-Beltrami operator on $(M,g).$
 
  For $k=1,\dots,n$ we define
\eq H_{k}=\fr{1}{\binom{n}{k}}\sum_{1\leq i_{1}<\dots<i_{k}\leq n}\k_{i_{1}}\cdots\k_{i_{k}}.\eeq
This includes the definition of the mean curvature,
\eq H=\fr 1n \sum_{i=1}^{n}\k_{i},\eeq
which deviates from some of the references. It corresponds to the notation in \cite{Roth2007}. Thus the traceless second fundamental form is
\eq \mathring{A}=A-Hg.\eeq
For smooth tensor fields on $M,$ $T=(t^{i_{1}\dots i_{k}}_{j_{1}\dots j_{l}}),$ we define the pointwise norms to be
\eq \|T\|=\sqrt{t^{i_{1}\dots i_{k}}_{j_{1}\dots j_{l}}t^{j_{1}\dots j_{l}}_{i_{1}\dots i_{k}}},\eeq 
where indices are lowered or lifted with respect to the induced metric of the hypersurface the tensor field is defined on. With the help of this definition we may define $L^{p}$-norms on a subset $\O\subset M$ to be
\eq \|T\|_{p,\O}=\(\int_{\O}\|T\|^{p}\)^{\fr 1p},\eeq
where the surface measure to be used is implicitly included in the set of integration $\O.$ Analogously we set
\eq \|T\|_{\infty,\O}=\sup_{\O}\|T\|.\eeq
The tensor $\mrm{Ric}=(R_{ij})$ is the Ricci tensor and $R=\mrm{tr}(\mrm{Ric})=R^{i}_{i}$ the scalar curvature. $\underline{\mrm{Ric}}(x)$ denotes the smallest eigenvalue of the Ricci tensor at $x\in M.$

For $M^{n}$ the symbol $\~{M}^{n}$ always denotes the normalized manifold
\eq \~{M}=|M|^{-\fr 1n}M\hra\R^{n+1}\eeq
with $|\~{M}|=1.$ The corresponding rescaled geometric quantities are denoted with a tilde as well, e.g.
\eq \~{g}=(\~{g}_{ij}),\ \~{A}=(\~{h}_{ij})\eeq
etc.

Finally
\eq B_{r}(x_{0})\subset\R^{n+1}\eeq
denotes an $(n+1)$-dimensional ball in $\R^{n+1}$ with radius $r$ and center $x_{0}.$ 

\section{Qualitative closeness revisited}\label{Qualitative closeness revisited}
In this section we turn our attention to the result, which connects $\l_{1}$ with closeness to a sphere, Theorem \ref{Roth}. We state, how the constant $C_{\e}$ involved here depends on $\e,$ whereafter we indicate, how this can be deduced from the corresponding sequence of lemmata in \cite{Roth2007}. We prove the following

\Prop \label{explepsilon}
In the situation of Theorem \ref{Roth} let $0<\e<\fr{2}{3\|H\|_{\infty}}.$ If (\ref{PC}) holds for 
\eq C_{\e}=\fr 12\min\(L\sqrt{\fr{n}{\l_{1}(M)}}\e^{2},L\),\eeq
where $L$ is bounded and uniformly positive whenever $\|H\|_{\infty}$ and $\|H_{2}\|_{2p}$ range in compact subsets of $(0,\infty),$ then we have
\eq M\subset B_{\sqrt{\fr{n}{\l_{1}(M)}}+\e}(x_{0})\bs B_{\sqrt{\fr{n}{\l_{1}(M)}}-\e}(x_{0}).\eeq
\eProp

\pf
We will spot and note the relevant formulae in \cite{Roth2007}, always denoting in which way they depend on the geometric quantities and on $\e.$ There is a sequence of constants, $C_{\e}$ is combined by. We start with \cite[p. 297, Lemma 2.1]{Roth2007}. First of all, it is required, that
\eq C_{\e}<\fr{n}{2}\|H_{2}\|_{2p}^{2}.\eeq Equation (5) yields
\eq A_{1}=\fr{2\|H\|_{\infty}^{2}}{\|H_{2}\|_{2p}^{2}}.\eeq
\cite[p. 298, Lemma 2.2]{Roth2007} yields
\eq A_{2}=\fr{A_{1}}{n\|H_{2}\|^{2}_{2p}}.\eeq
The proofs of \cite[Lemma 2.4, Lemma 2.5]{Roth2007} imply, that $A_{3}$ and $A_{4}$ are of a similar form. Finally, the author cites a lemma implying an $L^{\infty}$-estimate on the function
\eq \p=|X|\(|X|-\sqrt{\fr{n}{\l_{1}(M)}}\)^{2},\eeq
where $X$ is the position vector field with respect to the center of mass of $M,$ $x_{0}.$ The lemma is, cf. \cite[Lemma 3.1]{Roth2007},

\textit{For $p\geq 2$ and any $\eta>0,$ there exists $K_{\eta}(n,\|H\|_{\infty},\|H_{2}\|_{2p}),$ such that if (\ref{PC}) holds with $C_{\e}=K_{\eta},$ then $\|\p\|_{\infty}\leq \eta.$}

Essentially, the proof of this lemma is given in \cite[p. 16]{ColboisAndGrosjean09/2006}, also compare \cite[Sec. 6]{Roth2007}. Here one sees, that this $K_{\eta}$ can be chosen to be
\eq K_{\eta}=\min\(\fr{\eta}{(L'A_{4})^{4}}, c_{n}\)>0,\eeq
 where $L'$ is just of the same form as $A_{4}.$ Now, in \cite[p. 301]{Roth2007} the author defines
 \SAl \eta(\e)&=\min\(\(\sqrt{\fr{n}{\l_{1}(M)}}-\e\)\e^{2},\fr{1}{27\|H\|_{\infty}^{3}}\)\\
 			&\geq \min\(\fr 13 \sqrt{\fr{n}{\l_{1}(M)}}\e^{2},\fr{1}{27\|H\|_{\infty}^{3}}\),\end{split}\end{align}
 since $\e<\fr{2}{3\|H\|_{\infty}}$ and 
 \eq \l_{1}(M)\leq \fr{1}{n-1}\|R\|_{\infty}\leq n\|H\|_{\infty}^{2},\eeq
 compare \cite[Thm. 3.1]{Grosjean09/2002}. He concludes
 \eq M\subset B_{\sqrt{\fr{n}{\l_{1}(M)}}+\e}(x_{0})\bs B_{\sqrt{\fr{n}{\l_{1}(M)}}-\e}(x_{0})\eeq
under the assumption (\ref{PC}) with
\eq C_{\e}=\fr 12\min\(\fr n2\|H_{2}\|^{2}_{2p},c_{n},\fr{1}{3(L'A_{4})^{4}}\sqrt{\fr{n}{\l_{1}(M)}}\e^{2},\fr{1}{27(L'A_{4})^{4}\|H\|_{\infty}^{3}}\),\eeq
which has the form claimed in the proposition.
\epf

\section{Quantitative spherical closeness}\label{Quantitative spherical closeness}
Now we come to the proof of the main result. Let us state it again for a better readability.

\Thm\label{MainB}
Let $n\geq 2$ and $X\cn M^{n}\hra\R^{n+1}$ be the smooth, isometric embedding of a closed, connected, orientable and mean-convex hypersurface. Let $0<\a<1.$ Then there exists $c>0,$ such that whenever we have $\e<c|M|^{\fr{1}{n}}$ and the pointwise estimate 
\eq \label{PinchB} \|\mathring{A}\|\leq H|M|^{-\fr{2+a}{n}}\e^{2+\a}\eeq
holds, then $M$ is strictly convex and 
\eq \label{Maina} X(M)\subset B_{\sqrt{\fr{n}{\l_{1}(M)}}+\e}(x_{0})\bs B_{\sqrt{\fr{n}{\l_{1}(M)}}-\e}(x_{0}).\eeq
The constant $c$ depends on $n,$ $\a,$ $\|\~{A}\|_{\infty}$ and $\|\~{A}^{-1}\|_{\infty},$ where $|M|=\mrm{vol}(M),$ $\~{A}=|M|^{\fr{1}{n}}A,$ $\l_{1}(M)$ is the first nonzero eigenvalue of the Laplace-Beltrami operator on $M$ and $x_{0}$ is the center of mass of $M.$
\eThm

\pf
In this proof, $\~{C}_{i},$ $i\in\N,$ always denote generic constants which depend on $n,$ $\a,$ $\|\~{A}\|_{\infty}$ and $\|\~{A}^{-1}\|_{\infty}$ at most.
Set \eq p=n+1\eeq and let \eq k=\fr{6}{\a}.\eeq
For the rescaled surfaces 
\eq \~{M}=|M|^{-\fr 1n}M\eeq
we find from Theorem \ref{Perez}, that
\eq \label{Main3} \|\~{A}-\mu_{0}\~{g}\|_{kp}\leq \~{C}_{1}\|\mathring{\~{A}}\|_{kp},\eeq 
where $\mu_{0}=\mu_{0}(n,\a,\|\~{A}\|_{\infty},\|\~{A}^{-1}\|_{\infty}).$ 

The first condition we put on the constant $c$ is to satisfy
\eq c<\(\fr{1}{\sqrt{n(n-1)}}\)^{\fr{1}{2+\a}}.\eeq
Then (\ref{PinchB}) yields the strict convexity of $\~{M},$ due to \cite[Lemma 2.2]{Andrews:2012fk}. $\mu_{0}$ is strictly positive, since obviously we have
\eq \inf_{\~{M}}\~{\k}_{1}\leq \mu_{0}\leq \sup_{\~{M}}\~{\k}_{n}.\eeq

Define 
\eq \hat{M}=\mu_{0}\~{M}.\eeq
Then 
\eq \|\hat{A}-\hat{g}\|_{kp}=\(\int_{\hat{M}}\mu_{0}^{-kp}\|\~{A}-\mu_{0}\~{g}\|^{kp}\)^{\fr{1}{kp}}=\mu_{0}^{\fr{n}{kp}-1}\|\~{A}-\mu_{0}\~{g}\|_{kp}.\eeq
Define the set
\eq \hat{P}=\{\hat{x}\in\hat{M}\cn \|\hat{A}(\hat{x})-\hat{g}(\hat{x})\|<1\}.\eeq
Then its complement has volume
\eq |\hat{P}^{c}|\leq \int_{\hat{P}^{c}}\|\hat{A}-\hat{g}\|^{kp}\leq \mu_{0}^{n-kp}\|\~{A}-\mu_{0}\~{g}\|_{kp}^{kp}.\eeq
In order to apply Theorem \ref{Aubry}, we need an estimate on the Ricci tensor $\hat{\mrm{Ric}}=(\hat{R}_{ij}).$ By the Gaussian formula there holds
\eq \hat{R}_{ij}=n\hat{H}\hat{h}_{ij}-\hat{h}_{ik}\hat{h}^{k}_{j}.\eeq
Let $\hat{x}\in\hat{P}$ and $\xi\in T_{\hat{x}}\hat{M}.$ Then
\SAl \label{Main1} \hat{R}_{ij}\xi^{i}\xi^{j}&=n\hat{H}\hat{h}_{ij}\xi^{i}\xi^{j}-\hat{h}_{ik}\hat{h}^{k}_{j}\xi^{i}\xi^{j}\\
					&=n(\hat{H}-1)(\hat{h}_{ij}-\hat{g}_{ij})\xi^{i}\xi^{j}+n(\hat{h}_{ij}-\hat{g}_{ij})\xi^{i}\xi^{j}\\
					&\hp{=} +n(\hat{H}-1)\|\xi\|^{2}+(n-1)\|\xi\|^{2}-2(\hat{h}_{ij}-\hat{g}_{ij})\xi^{i}\xi^{j}\\
					&\hp{=}-(\hat{h}_{ik}-\hat{g}_{ik})(\hat{h}^{k}_{j}-\d^{k}_{j})\xi^{i}\xi^{j},\end{split}\end{align}
from which we obtain at $\hat{x}$
\eq \|\hat{\mrm{Ric}}-(n-1)\hat{g}\|\leq \~{C}_{2}\|\hat{A}-\hat{g}\|,\eeq
since $\|\hat{A}-\hat{g}\|<1.$ In the notation of Theorem \ref{Aubry} we obtain
\SAl \int_{\hat{M}}(\underline{\hat{\mrm{Ric}}}-(n-1))^{kp}_{-}&\leq \int_{\hat{P}}\~{C}_{2}^{kp}\|\hat{A}-\hat{g}\|^{kp}+\int_{\hat{P}^{c}}(\underline{\hat{\mrm{Ric}}}-(n-1))_{-}^{kp}\\
									&\leq \(\~{C}_{2}^{kp}\mu_{0}^{n-kp}+(n-1)^{kp}\mu_{0}^{n-kp}\)\|\~{A}-\mu_{0}\~{g}\|_{kp}^{kp}\\
									&= \~{C}_{3}\|\~{A}-\mu_{0}\~{g}\|_{kp}^{kp}.\end{split}\end{align} 
Thus Theorem \ref{Aubry} will be applicable under condition (\ref{PinchB}), if we choose $c$ small enough to ensure the last of the following inequalities (note that in the first inequality we use (\ref{Main3})).
\SAl \~{C}_{3}\|\~{A}-\mu_{0}\~{g}\|^{kp}_{kp}&\leq \~{C}_{3}\~{C}_{1}^{kp}\|\mathring{\~{A}}\|^{kp}_{kp}=\~{C}_{3}\~{C}_{1}^{kp}|M|^{\fr{kp}{n}-1}\|\mathring{A}\|^{kp}_{kp}\\
								&\leq \~{C}_{3}\~{C}_{1}^{kp}|M|^{-\fr{(1+\a)kp+n}{n}}\e^{(2+\a)kp}\|H\|^{kp}_{kp}\\
								&= \~{C}_{3}\~{C}_{1}^{kp}|M|^{-\fr{(2+\a)kp}{n}}\e^{(2+\a)kp}\|\~{H}\|^{kp}_{kp}\\
								&<\~{C}_{3}\~{C}_{1}^{kp}c^{(2+\a)kp}\|\~{H}\|_{kp}^{kp}\\
								&\overset{!}{<}\fr{|\hat{M}|}{C(n,kp)}=\fr{\mu_{0}^{n}}{C(n,kp)},\end{split}\end{align}
where $C(n,kp)$ is the constant from Theorem \ref{Aubry}. Thus $c=c(n,\a,\|\~{A}\|_{\infty},\|\~{A}^{-1}\|_{\infty})$ is additionally choosable, such that this chain of inequalities is true. We may apply Theorem \ref{Aubry} to conclude
\SAl \l_{1}(\hat{M})&\geq n\(1-C(n,kp)\(\fr{1}{|\hat{M}|}\int_{\hat{M}}(\underline{\hat{\mrm{Ric}}})-(n-1))_{-}^{kp}\)^{\fr{1}{kp}}\)\\
			&\geq n\(1-C(n,kp)\mu_{0}^{-\fr{n}{kp}}\~{C}_{1}\~{C}_{3}^{\fr{1}{kp}}\|\~{H}\|_{kp}\~{\e}^{2+\a}\),\end{split}\end{align}
where $\~{\e}=|M|^{-\fr 1n}\e.$ We obtain
\eq \label{Main4} \l_{1}(\~{M})\geq \mu_{0}^{2}n(1-\~{C}_{4}\~{\e}^{2+\a}),\eeq
with a new constant $\~{C}_{4}.$

Now we want to apply Theorem \ref{Roth}. Therefore we need estimates of the curvature integrals. First note, that
\eq \~{H}_{2}=\fr{1}{n(n-1)}\~{R}.\eeq
A similar calculation as (\ref{Main1}) shows, that at any point 
\eq \~{x}\in\~{P}_{\g}=\{\~{x}\in\~{M}\cn\|\~{A}-\mu_{0}\~{g}\|<\g\},\ 0<\g<1,\eeq
we have
\eq \label{Main5} \|\~{R}_{ij}-\mu_{0}^{2}(n-1)\~{g}_{ij}\|\leq \~{C}_{5}(n,\mu_{0})\|\~{A}-\mu_{0}\~{g}\|.\eeq 
Furthermore there holds
\eq |\~{P}_{\g}^{c}|\g^{kp}\leq \int_{\~{P}_{\g}^{c}}\|\~{A}-\mu_{0}\~{g}\|^{kp}\leq \~{C}_{1}^{kp}\|\mathring{\~{A}}\|^{kp}_{kp}\leq \~{C}_{6}\~{\e}^{(2+\a)kp}\eeq
and thus
\eq |\~{P}_{\g}^{c}|\leq \~{C}_{6}\(\fr{\~{\e}^{2+\a}}{\g}\)^{kp}.\eeq
We estimate 
\SAl\label{Main8} \(\int_{\~{M}}\~{H}^{2p}_{2}\)^{\fr 1p}&=\(\int_{\~{P}_{\g}}\(\fr{\~{R}}{n(n-1)}\)^{2p}+\int_{\~{P}^{c}_{\g}}\(\fr{\~{R}}{n(n-1)}\)^{2p}\)^{\fr 1p}\\
							&\leq \left\|\fr{\~{R}}{n(n-1)}\right\|^{2}_{2p,\~{P}_{\g}}+\left\|\fr{\~{R}}{n(n-1)}\right\|^{2}_{2p,\~{P}^{c}_{\g}}\\
							&\leq \(\mu_{0}^{2}+\~{C}_{5}\|\~{A}-\mu_{0}\~{g}\|_{2p,\~{P}_{\g}}\)^{2}+|\~{P}^{c}_{\g}|^{\fr 1p}\|\~{H}\|^{4}_{\infty},\end{split}\end{align}
where we used $\~{H}_{2}^{\fr 12}\leq \~{H}$ and (\ref{Main5}).

Furthermore we obtain from (\ref{Main5}), that
\SAl\label{Main6} \(\int_{\~{M}}\~{H}\)^{2}\geq \(\int_{\~{P}_{\g}}\(\fr{\~{R}}{n(n-1)}\)^{\fr 12}\)^{2}&\geq \(|\~{P}_{\g}|\sqrt{\mu_{0}^{2}-\~{C}_{5}\g}\)^{2}\\
					&=|\~{P}_{\g}|^{2}\mu_{0}^{2}-|\~{P}_{\g}|^{2}\~{C}_{5}\g \end{split}\end{align}
for all
\eq 0<\g<\fr{\mu_{0}^{2}}{\~{C}_{5}}.\eeq
From (\ref{Main4}), (\ref{Main8}) and (\ref{Main6}) we obtain
\SAl &\l_{1}(\~{M})\(\int_{\~{M}}\~{H}\)^{2}-n\|\~{H}_{2}\|_{2p}^{2}\\
	\geq &(\mu_{0}^{2}n-\mu_{0}^{2}n\~{C}_{4}\~{\e}^{2+\a})(|\~{P}_{\g}|^{2}\mu_{0}^{2}-|\~{P}_{\g}|^{2}\~{C}_{5}\g)\\
		&\hp{=}-n\mu_{0}^{4}-n\~{C}_{5}^{2}\g^{2}-2n\mu_{0}^{2}\~{C}_{5}\g-n|\~{P}^{c}_{\g}|^{\fr 1p}\|\~{H}\|^{4}_{\infty}\\
		\geq &-\~{C}_{7}|\~{P}^{c}_{\g}|-\~{C}_{7}\g-\~{C}_{7}\~{\e}^{2+\a}-\~{C}_{7}\(\fr{\~{\e}^{2+\a}}{\g}\)^{k},\end{split}\end{align}
where $\~{C}_{7}$ is a new constant. According to Theorem \ref{Roth} and Proposition \ref{explepsilon} there exists $C_{\~{\e}},$ which can be chosen as 
\eq C_{\~{\e}}=\fr 12\min\(L\sqrt{\fr{n}{\l_{1}(\~{M})}}\~{\e}^{2},L\),\eeq
 such that whenever $\~{\e}<\fr{2}{3\|\~{H}\|_{\infty}}$ and 
 \eq \l_{1}(\~{M})\(\int_{\~{M}}\~{H}\)^{2}-n\|\~{H}_{2}\|^{2}_{2p}>-C_{\~{\e}},\eeq
 we could conclude 
 \eq \~{M}\subset B_{\sqrt{\fr{n}{\l_{1}(\~{M})}}+\~{\e}}(\~{x}_{0})\bs B_{\sqrt{\fr{n}{\l_{1}(\~{M})}}-\~{\e}}(\~{x}_{0}).\eeq
  Now define
 \eq \g=\~{\e}^{2+\fr{\a}{2}}.\eeq
 Then
 \SAl &\~{C}_{7}\(\(\fr{\~{\e}^{2+\a}}{\g}\)^{kp}+\(\fr{\~{\e}^{2+\a}}{\g}\)^{k}+\g+\~{\e}^{2+\a}\)\\
 	\leq &\~{C}_{7}\(\~{\e}^{\fr{\a kp}{2}}+\~{\e}^{\fr{\a k}{2}}+\~{\e}^{2+\fr{\a}{2}}+\~{\e}^{2+\a}\)\\
 												=&\~{C}_{7}\(\~{\e}^{3p}+\~{\e}^{3}+\~{\e}^{2+\fr{\a}{2}}+\~{\e}^{2+\a}\)\\
												<&\fr{1}{2}\min\(L\sqrt{\fr{n}{\l_{1}(\~{M})}}\~{\e}^{2},L\),\end{split}\end{align}
for all $0<\~{\e}<c,$ if $c$ is small enough in dependence of $n,$ $\a,$ $\|\~{A}\|_{\infty}$ and $\|\~{A}^{-1}\|_{\infty},$ such that the requirements for $\g,$ namely
\eq \g<\min\(1,\fr{\mu_{0}^{2}}{\~{C}_{5}}\),\eeq
are fulfilled as well. 

We conclude, rescaling again,
\eq M\subset B_{\sqrt{\fr{n}{\l_{1}(M)}}+\e}(x_{0})\bs B_{\sqrt{\fr{n}{\l_{1}(M)}}-\e}(x_{0}),\eeq
the desired result.
\epf

\Rem
The previous result is easier to comprehend, if one restricts to the class of hypersurfaces of bounded volume and modulus of convexity, namely
\eq 0<c\leq |M|\leq C\eeq
and 
\eq 0<cg\leq A\leq Cg.\eeq
Then, in order to prove $\e$-closeness, one has to find constants $c>0$ and $\b>0,$ such that
\eq \|A-Hg\|\leq cH\e^{2+\b},\eeq
where $c$ must not depend on $\e.$ Then applying Theorem \ref{MainB} with $\a=\fr{\b}{2},$ one concludes $\e$-closeness for small $0<\e<\e_{0}.$
\eRem

\section{Concluding remarks and open questions}
We must not hesitate to remark, that this result is only a first step towards a better understanding of the stability problem. It helps to control the order of the sufficient $\d$ with respect to $\e,$ which is sufficient for first applications in geometric flows, compare \cite{js:Ipcf}. 

However, two things will be desirable in this context. Firstly, there would be direct applications to geometric flows, if one could improve the order $\e^{2+\a}.$ We are not aware of the existence of such a result. Secondly, pinching results for the first eigenvalue of the Laplacian are known in other ambient space, cf. \cite{GrosjeanAndRoth2012}. It would be interesting, with immediate applications to curvature flows in those spaces, whether results like ours could be deduced in those settings as well.

\bibliographystyle{hamsplain}
\bibliography{MyBiblio}

\end{document}